\documentclass[12pt]{article}
\usepackage{graphicx}
\usepackage{amsmath,amsthm,amssymb,enumerate}
\usepackage{euscript,mathrsfs}
\usepackage{color}
\usepackage{dsfont}
\usepackage[left=2cm,right=2cm,top=3.5cm,bottom=3.5cm]{geometry}
\usepackage{color}
\usepackage[framemethod=tikz]{mdframed}
\allowdisplaybreaks

\usepackage{soul}

\catcode`\@=11 \@addtoreset{equation}{section}

\catcode`\@=12

\allowdisplaybreaks

\newtheorem{Theorem}{Theorem}[section]
\newtheorem{Proposition}[Theorem]{Proposition}
\newtheorem{Lemma}[Theorem]{Lemma}
\newtheorem{Corollary}[Theorem]{Corollary}

\theoremstyle{definition}
\newtheorem{Definition}[Theorem]{Definition}

\newtheorem{Remark}[Theorem]{Remark}

\newcommand{\bTheorem}[1]{
\begin{Theorem} \label{T#1} }
\newcommand{\eT}{\end{Theorem}}

\newcommand{\bProposition}[1]{
\begin{Proposition} \label{P#1}}
\newcommand{\eP}{\end{Proposition}}

\newcommand{\bLemma}[1]{
\begin{Lemma} \label{L#1} }
\newcommand{\eL}{\end{Lemma}}

\newcommand{\bCorollary}[1]{
\begin{Corollary} \label{C#1} }
\newcommand{\eC}{\end{Corollary}}

\newcommand{\bRemark}[1]{
\begin{Remark} \label{R#1} }
\newcommand{\eR}{\end{Remark}}

\newcommand{\bDefinition}[1]{
\begin{Definition} \label{D#1} }
\newcommand{\eD}{\end{Definition}}

\newcommand{\Del}{\Delta_x}

\newcommand{\vrB}{\vr_B}

\newcommand{\Ds}{\mathbb{D}_x}

\newcommand{\vuB}{\vc{u}_B}

\newcommand{\bfphi}{\boldsymbol{\varphi}}

\newcommand{\bFormula}[1]{
\begin{equation} \label{#1}}
\newcommand{\eF}{\end{equation}}

\newcommand{\Ov}[1]{\overline{#1}}

\newcommand{\aleq}{\stackrel{<}{\sim}}

\newcommand{\vr}{\varrho}

\newcommand{\tvt}{\tilde \vt}

\newcommand{\vt}{\vartheta}
\newcommand{\vu}{\vc{u}}
\newcommand{\vm}{\vc{m}}

\newcommand{\vc}[1]{{\bf #1}}

\newcommand{\Div}{{\rm div}_x}
\newcommand{\Grad}{\nabla_x}

\newcommand{\dx}{\,{\rm d} {x}}

\newcommand{\dt}{\,{\rm d} t }

\newcommand{\intO}[1]{\int_{\Omega} #1 \ \dx}

\newcommand{\intgo}[1]{\int_{\partial \Omega} #1 \ [ \vuB \cdot \vc{n} ]^+ \D \sigma_x}
\newcommand{\intgi}[1]{\int_{\partial \Omega} #1 \ [ \vuB \cdot \vc{n} ]^- \D \sigma_x}

\newcommand{\D}{{\rm d}}

\newcommand{\ep}{\varepsilon}

\newcommand{\R}{\mathbb{R}}
\newcommand{\vtB}{\vt_B}

\newcommand{\br}{ \nonumber \\ }

\def\softd{{\leavevmode\setbox1=\hbox{d}%
          \hbox to 1.05\wd1{d\kern-0.4ex{\char039}\hss}}}
\definecolor{Cgrey}{rgb}{0.85,0.85,0.85}
\definecolor{Cblue}{rgb}{0.50,0.85,0.85}
\definecolor{Cred}{rgb}{1,0,0}
\definecolor{fancy}{rgb}{0.10,0.85,0.10}

\newcommand\Cbox[2]{%
    \newbox\contentbox%
    \newbox\bkgdbox%
    \setbox\contentbox\hbox to \hsize{%
        \vtop{
            \kern\columnsep
            \hbox to \hsize{%
                \kern\columnsep%
                \advance\hsize by -2\columnsep%
                \setlength{\textwidth}{\hsize}%
                \vbox{
                    \parskip=\baselineskip
                    \parindent=0bp
                    #2
                }%
                \kern\columnsep%
            }%
            \kern\columnsep%
        }%
    }%
    \setbox\bkgdbox\vbox{
        \color{#1}
        \hrule width  \wd\contentbox %
               height \ht\contentbox %
               depth  \dp\contentbox
        \color{black}
    }%
    \wd\bkgdbox=0bp%
    \vbox{\hbox to \hsize{\box\bkgdbox\box\contentbox}}%
    \vskip\baselineskip%
}

\mdfdefinestyle{MyFrame}{%
	linecolor=black,
	outerlinewidth=1pt,
	roundcorner=5pt,
	innertopmargin=\baselineskip,
	innerbottommargin=\baselineskip,
	innerrightmargin=10pt,
	innerleftmargin=10pt,
	backgroundcolor=white!20!white}

\date{}

\begin{document}


\title{Asymptotic stability of solutions to the Navier--Stokes--Fourier system driven by inhomogeneous Dirichlet boundary conditions}

\author{Eduard Feireisl
\thanks{The work of E. Feireisl was partially supported by the
Czech Sciences Foundation (GA\v CR), Grant Agreement
21--02411S. The Institute of Mathematics of the Academy of Sciences of
the Czech Republic is supported by RVO:67985840.} \and Young--Sam Kwon
\thanks{The work of Y.--S. Kwon was partially supported by
the National Research Foundation of Korea (NRF2020R1F1A1A01049805)}}

\date{\today}

\maketitle

\bigskip

\centerline{Institute of Mathematics of the Academy of Sciences of the Czech Republic}

\centerline{\v Zitn\' a 25, CZ-115 67 Praha 1, Czech Republic}

\medskip

\centerline{Department of Mathematics, Dong-A University}

\centerline{Busan 49315, Republic of Korea}

\bigskip

\begin{abstract}
	
	We consider global in time solutions of the Navier--Stokes--Fourier system describing the motion of a general compressible, viscous and heat conducting fluid far from equilibirum. Using a new concept of weak solution suitable to accommodate the inhomogeneous Dirichlet time dependent data we find sufficient conditions for the global in time weak solutions to be ultimately bounded.

\end{abstract}

{\bf Keywords:} Navier--Stokes--Fourier system, long--time behavior, bounded absorbing set, Dirichlet boundary conditions

\bigskip


\section{Introduction}
\label{i}

The Navier--Stokes--Fourier system describing the time evolution of the mass density $\vr = \vr(t,x)$, the velocity $\vu = \vu(t,x)$, and the temperature $\vt = \vt(t,x)$ of a general compressible, viscous, and heat--conducting fluid, endowed with inhomogeneous boundary conditions, is a prominent example of a \emph{dissipative system} in the framework of continuum fluid mechanics. In general, a dissipative system is a thermodynamically open system
confined to a physical space $\Omega \subset R^d$ and
considered far from equilibrium, exchanging energy and matter with the outer world. The field equations of the \emph{Navier--Stokes--Fourier system} describing the motion in the interior of the cavity $\Omega$ read:
\begin{align}
\partial_t \vr + \Div (\vr \vu) &= 0, \label{i1} \\
	\partial_t (\vr \vu) + \Div (\vr \vu \otimes \vu) + \Grad p(\vr, \vt) &= \Div \mathbb{S} (\vt, \Ds \vu) + \vr \vc{g},
	\label{i2} \\
\partial_t (\vr e(\vr, \vt)) + \Div (\vr e(\vr, \vt) \vu ) + \Div \vc{q}(\vt, \Grad \vt) &= \mathbb{S} (\vt, \Ds \vu) : \Ds \vu - p(\vr, \vt) \Div \vu,
\label{i3} 	
	\end{align}
where the \emph{viscous stress} $\mathbb{S}$ is given by Newton's law
\begin{equation} \label{i4}
	\mathbb{S} (\vt, \Ds \vu) = \mu(\vt) \left( \Grad \vu + \Grad^t \vu - \frac{2}{d} \Div \vu \mathbb{I} \right) + \eta (\vt) \Div \vu \mathbb{I}, \ \Ds \vu \equiv \frac{1}{2} (\Grad \vu + \Grad^t \vu),
	\end{equation}
and the \emph{heat flux} $\vc{q}$ by Fourier's law
\begin{equation} \label{i5}
	\vc{q} (\vt, \Ds \vu) = - \kappa (\vt) \Grad \vt .
\end{equation}
The \emph{pressure} $p$ and the \emph{internal energy} $e$ are interrelated through Gibbs' equation
\begin{equation} \label{i6}
	\vt D s = D e + p D \left( \frac{1}{\vr} \right),
\end{equation}
where $s$ is the \emph{entropy}. In view of \eqref{i6}, the internal energy balance \eqref{i3} may be replaced by the entropy balance
\begin{equation} \label{i7}
\partial_t (\vr s(\vr, \vt)) + \Div (\vr s(\vr, \vt) \vu ) + \Div \left( \frac{\vc{q}(\vt, \Grad \vt)}{\vt} \right) = \frac{1}{\vt} \left( \mathbb{S} (\vt, \Ds \vu) : \Ds \vu - \frac{ \vc{q} \cdot \Grad \vt }{\vt} \right).
	\end{equation}
As we shall see below, it is the entropy balance \eqref{i7} that is more convenient for the \emph{weak} formulation of the problem.

The fluid occupies a bounded domain $\Omega$, whereas the velocity $\vu$ and the temperature $\vt$ satisfy the inhomogeneous Dirichlet boundary conditions
\begin{align}
	\vu &= \vuB \ \mbox{on} \ (T,\infty) \times \partial \Omega, \label{i8} \\
	\vt &= \vtB \ \mbox{on} \ (T,\infty) \times \partial \Omega, \label{i9}.
	\end{align}
 Accordingly, the boundary of the space--time cylinder $(T,\infty) \times \partial \Omega$ can be decomposed as
 \begin{align}
 	\Gamma_{\rm in} & =\left\{ (t,x) \in (T,\infty) \times \partial \Omega \ \Big|\ \vuB(t,x) \cdot \vc{n}(x) < 0 \right\}, \br
\Gamma_{\rm wall} & =\left\{ (t,x) \in (T,\infty) \times \partial \Omega \ \Big|\ \vuB(t,x)\cdot \vc{n}(x) = 0 \right\}, \br
\Gamma_{\rm out} & = \left\{ (t,x) \in (T,\infty) \times \partial \Omega \ \Big|\ \vuB(t,x) \cdot \vc{n}(x) > 0 \right\},	
\nonumber
\end{align}
where $\vc{n}$ denotes the outer normal vector to $\partial \Omega$. Finally, the density must be prescribed on the inflow component of the boundary,
\begin{equation} \label{i10}
	\vr = \vrB \ \mbox{on}\ \Gamma_{\rm in}.
	\end{equation}

As shown in a series of papers by Matsumura and Nishida \cite{MANI1}, \cite{MANI}, Valli \cite{Vall2}, \cite{Vall1}, Valli and Zajaczkowski \cite{VAZA}, the Navier--Stokes--Fourier system is globally well posed
in the class of classical solutions
in a small neighbourhood of a stable equilibrium. Unfortunately, this perturbation technique cannot be used in the \emph{far from equilibrium regime}
usually associated to turbulence.
In view of the well known and so far unsurmountable problems concerning the existence
of suitable {\it a priori} bounds for large data global in time solutions to
non--linear problems in fluid dynamics, the only available framework
are the weak solutions in the spirit of the pioneering work by Leray \cite{LER}, and later, in the context of compressible fluids, by Lions \cite{LI4}.

A mathematical theory based on weak solutions for the complete fluid systems was presented in
\cite{FeNo6A}. Unfortunately, the concept of weak solutions developed in \cite{FeNo6A} applies to
\emph{energetically closed} systems with $\vuB = 0$ and with \eqref{i9} replaced by the homogeneous Neumann boundary condition
\[
\Grad \vt \cdot \vc{n}|_{\partial \Omega} = 0.
\]
Under these circumstances, the long time behaviour of solutions is well understood, and, in the case of a time independent driving force $\vc{g} = \vc{g}(x)$, obeys the following dichotomy:
\begin{itemize}
\item
Either $\vc{g} = \Grad F$ and then all solutions tend an equilibrium;
\item
or  $\vc{g} \ne \Grad F$ and then
\[
\intO{ \left[ \frac{1}{2} \vr |\vu|^2 + \vr e(\vr, \vt) \right](t, \cdot) } \to \infty \ \mbox{as}\ t \to \infty
\]
for any weak solution $(\vr, \vu, \vt)$;

\end{itemize}
see \cite{FP20}, \cite{FeiPr}. It turns out that driving the system by means of a non--trivial volume force while keeping the conservative boundary conditions is not realistic and definitely not suitable for describing
phenomena related to turbulence.

The mathematical theory of weak solutions has been extended to the \emph{energetically open} Navier--Stokes--Fourier system only recently in \cite{FeiNov20}, and, finally, in \cite{ChauFei}. In particular, the inhomogeneous boundary condition for
the temperature requires a new approach, developed in \cite{ChauFei}, based on the balance of ballistic energy in the sense of Ericksen \cite{Eri}. Recently, the existence of of weak solutions has also been proved for  a bi-fluid model for a mixture of two compressible
non interacting fluids with general boundary data in \cite{KKNN}.

 Note that this kind of boundary conditions is physically relevant, in particular
for the Rayleigh--B\' enard problem and Taylor experiment arising in models of turbulence, cf. Birnir and Svanstedt \cite{BirSva}, Constantin et al. \cite{CFT}, \cite{CFNT}, Davidson \cite{DAVI} among others.

To the best of our knowledge, this is the first attempt to describe the asymptotic behaviour of the Navier--Stokes--Fourier system with the \emph{inhomogeneous Dirichlet} boundary conditions. We focus on the problem of \emph{global boundedness} of trajectories and the existence of bounded absorbing sets.
This is the concept of \emph{dissipativity in the sense of Levinson} extending the classical approach to closed systems via a Lyapunov function,
see e.g. Haraux \cite{Hara1}, Kuznetsov and Reitmann \cite[Chapter 1, Section 1.2]{KuzReit}. The crucial quantity
if the \emph{ballistic energy}
\[
\intO{ E_{B} (\vr, \vt, \vu) } ,\ E_B (\vr, \vt, \vu) \equiv  \frac{1}{2} \vr |\vu - \vuB|^2 + \vr e (\vr, \vt) - \vtB \vr s(\vr, \vt)
\]
where $\vuB$, $\vtB$ are suitable extensions of the boundary data inside $\Omega$.

\medskip

{\it We say that the Navier--Stokes--Fourier system \eqref{i1}--\eqref{i10} is Levinson dissipative if there exists a universal constant $\mathcal{E}_\infty$ such that}
\[
\limsup_{t \to \infty} \intO{ E_B (\vr, \vt, \vu) } \leq \mathcal{E}_\infty
\]
{\it for any (weak) solution $(\vr, \vt, \vu)$ defined on a time interval $(T_0, \infty)$.}

\medskip

Our goal is to identify a class of constitutive relations (equations of state (EOS), viscosity and heat conductivity coefficients), for which the Navier--Stokes--Fourier system \eqref{i1}--\eqref{i10}
is Levinson dissipative. We adopt the following strategy:

\begin{enumerate}
	\item In Section \ref{w}, we recall the concept of weak solution introduced in \cite{ChauFei}, together with the associated ballistic energy inequality.
	\item Inspired by \cite{CiFeJaPe1}, we introduce a class of equations of state penalizing the pressure if the density approaches a critical value $\Ov{\vr}$, see Section \ref{h}.
	\item The main results are stated in Section \ref{m}.
	\item In Section \ref{u}, we show the uniform bounds of the ballistic energy proving Levinson dissipativity of the Navier--Stokes--Fourier system. We also study convergence to
	equilibrium solutions for a particular class of boundary data.
	
	\item Possible applications including the existence of global attractors, statistical solutions, and the existence of time periodic solutions are briefly sketched in Section \ref{a}.
	\end{enumerate}

\section{Weak solutions, ballistic energy}
\label{w}

We recall the concept of weak solution and ballistic energy introduced in \cite{ChauFei}.

\subsection{Weak solutions}

We consider the Navier--Stokes--Fourier system \eqref{i1}--\eqref{i10}
defined on the set $(T,\infty) \times \Omega$, where $T < \infty$. Although the existence theory
is formulated in terms of the initial data
\[
	\vr(T, \cdot) = \vr_0, \ (\vr \vu)(T, \cdot) = \vm_0,\
	\vr s (T, \cdot) = S_0,\ S_0 = \vr_0 s(\vr_0, \vt_0),
	\]
their specific form is irrelevant for the analysis of the present paper.

\begin{mdframed}[style=MyFrame]
	
	\begin{Definition}[Weak solution] \label{Dw1}
		
		\medskip
		We say that a trio of functions $(\vr, \vt, \vu)$ is a \emph{weak solution} of the Navier--Stokes--Fourier system \eqref{i1}--\eqref{i10} in $(T, \infty) \times \Omega$
		if the following holds:
		\begin{itemize}

			\item {\bf Equation of continuity.}
			$\vr \geq 0$ and the integral identity
			\begin{align}
				\int_T^\infty &\intO{ \Big[ \vr \partial_t \varphi + \vr \vu \cdot \Grad \varphi \Big]} \dt \br
				&= \int_T^\infty \intgi{ \varphi \vrB } \dt + \int_T^\infty \intgo{ \varphi \vr} \dt
				\label{w2}
			\end{align}
				
		\noindent
			holds for any $\varphi \in C^1_c((T,\infty) \times \Ov{\Omega})$.
			
			\item {\bf Momentum equation.}
			\begin{align}
				\vu &\in L^r_{\rm loc}(T,\infty; W^{1,r}(\Omega; R^d)) \ \mbox{for some}\ r > 1,\br (\vu - \vuB) &\in L^r_{\rm loc}(T,\infty; W^{1,r}_0(\Omega; R^d)),
				\label{w3}		
			\end{align}
			and	
				\begin{align}
					\int_T^\infty &\intO{ \Big[ \vr \vu \cdot \partial_t \bfphi + \vr \vu \otimes \vu : \Grad \bfphi + p(\vr, \vt) \Div \bfphi \Big] } \dt \br &=
					\int_T^\infty \intO{ \Big[ \mathbb{S} (\vt, \Ds \vu) : \Ds \bfphi - \vr \vc{g} \cdot \bfphi \Big] } \dt
					\label{w4}
				\end{align}
				for any $\bfphi \in C^1_c((T,\infty) \times \Omega; R^d)$.
			
			\item {\bf Entropy inequality.}	
				\begin{align}
					&  - \int_{T}^{\infty} \intO{
						\left[ \vr s \partial_t \varphi + \vr s \vu \cdot \Grad \varphi + \frac{\vc{q}}{\vt} \cdot \Grad \varphi \right] } \dt \br &\geq
					\int_{T}^{\infty} \intO{ \frac{\varphi}{\vt} \left( \mathbb{S}(\vt, \Ds \vu ): \Ds \vu - \frac{\vc{q}(\vt, \Grad \vt) \cdot \Grad \vt }{\vt} \right) } \dt
					\label{w5}
				\end{align}
				
			\noindent
			for any $\varphi \in C^1_c((T,\infty) \times \Omega)$, $\varphi \geq 0$.
			
			\item {\bf Ballistic energy inequality.}
			For any
			\[
			\tvt \in C^1([T,\infty) \times \Ov{\Omega}),\ \tvt > 0,\ \tvt|_{\partial \Omega} = \vtB
			\]
			there holds 		
				\begin{align}
				&-
				\int_T^\infty \partial_t \psi \intO{ \left( \frac{1}{2} \vr |\vu - \vuB|^2 + \vr e - \tvt \vr s \right) } \dt  \br &+
				\int_T^\infty \psi \intgi{ \Big[ \vrB e(\vrB, \vtB) - \vtB \vrB s(\vrB, \vtB) \Big] } \dt \br &+
				\int_T^\infty \psi \intgo{ \Big[ \vr e(\vr, \vtB) - \vtB \vr s(\vr, \vtB) \Big] } \dt \br &+
				\int_T^\infty \psi \intO{ \frac{\tvt}{\vt}	 \left( \mathbb{S} : \Ds \vu - \frac{\vc{q} \cdot \Grad \vt }{\vt} \right) } \dt \br
				&\leq
				- \int_T^\infty \psi \intO{ \Big[ \vr (\vu - \vuB) \otimes (\vu - \vuB) + p \mathbb{I} - \mathbb{S} \Big] : \Ds \vuB } \dt \br &+
				\int_T^\infty \psi \intO{ \vr (\vu - \vuB)\cdot (\vc{g} - \partial_t \vuB - \vuB \cdot \Grad \vuB) } \dt \br
				&- \int_T^\infty \psi \intO{ \left[ \vr s \left( \partial_t \tvt + \vu \cdot \Grad \tvt \right) + \frac{\vc{q}}{\vt} \cdot \Grad \tvt \right] } \dt
				\label{w6}
			\end{align}
			for any $\psi \in C^1_c (T,\infty)$, $\psi \geq 0$.

		\end{itemize}	
		
	\end{Definition}

	\end{mdframed}

The symbol $\vuB$ in \eqref{w6} denotes any $C^1$ extension of the boundary velocity $\vuB$. It can be shown (cf. \cite[Remark 2.2]{BreFeiNov20}) that the specific form of the ballistic energy
inequality is independent of the extension $\vuB$. Specifically, if \eqref{w6} holds for some $\vuB$, then it holds for any $\vuB$ attaining the same boundary value. We point out that
such a statement may not hold for the temperature extension $\tvt$.

\begin{Remark} \label{Rw1}
	
	The regularity of a weak solution $(\vr, \vt, \vu)$ is determined by the available {\it a priori} bounds based mostly on the ballistic energy inequality. In particular, the density $\vr$ is only Lebesgue integrable while both
	\eqref{w2} and \eqref{w6} refer to its trace on $\Gamma_{\rm out}$. The latter is understood in the following way. The velocity $\vu$ is a Sobolev function that admits a trace. Moreover, the vector field
	$[\vr, \vr \vu]$ having zero space--time divergence admits a normal trace on the space time cylinder $(T, \infty) \times \Omega$. Consequently, $\vr|_{\Gamma_{\rm out}}$ is determined through
	\[
	\vr \vu \cdot \vc{n} |_{\partial \Omega} = \vr|_{\Gamma_{\rm out}} \vuB \cdot \vc{n}.
		\]	
	 The reader may consult \cite{FeiNov20} for details.
	
	\end{Remark}

The existence of global--in--time weak solutions under certain restrictions imposed on the constitutive relations was proved in \cite[Theorem 4.2]{ChauFei}. The weak solutions also comply with the weak--strong uniqueness principle, see \cite[Theorem 3.1]{ChauFei}. Specifically, any weak solution in the sense of Definition \ref{Dw1} coincides with the strong solution of the Navier--Stokes--Fourier system driven by the same initial/boundary data as long as the strong solution exists.

The ballistic energy inequality \eqref{w6} can be written in a more concise form
	\begin{align}
	\frac{\D }{\dt} &\intO{ \left( \frac{1}{2} \vr |\vu - \vuB|^2 + \vr e - \tvt \vr s \right) }   \br &+
	\intgi{ \Big[ \vrB e(\vrB, \vtB) - \vtB \vrB s(\vrB, \vtB) \Big] }  \br &+
	\intgo{ \Big[ \vr e(\vr, \vtB) - \vtB \vr s(\vr, \vtB) \Big] } \br &+
	\intO{ \frac{\tvt}{\vt}	 \left( \mathbb{S} : \Ds \vu - \frac{\vc{q} \cdot \Grad \vt }{\vt} \right) }  \br
	&\leq - \intO{ \Big[ \vr (\vu - \vuB) \otimes (\vu - \vuB) + p \mathbb{I} - \mathbb{S} \Big] : \Ds \vuB }  \br &+
	\intO{ \vr (\vu - \vuB)\cdot (\vc{g} - \partial_t \vuB - \vuB \cdot \Grad \vuB) }  \br
	&-  \intO{ \left[ \vr s \left( \partial_t \tvt + \vu \cdot \Grad \tvt \right) + \frac{\vc{q}}{\vt} \cdot \Grad \tvt \right] }
	\label{w7}
\end{align}
understood in $\mathcal{D}'(T, \infty)$.

Our goal is to identify the class of boundary data $\vuB$, $\vtB$, for which \eqref{w7} gives rise to a globally bounded ballistic energy,
\begin{equation} \label{w8}
\limsup_{t \to \infty} \intO{ \left( \frac{1}{2} \vr |\vu - \vuB|^2 + \vr e - \vtB \vr s \right) } \leq \mathcal{E}_\infty.
\end{equation}
for a certain $\vtB$. In other words, the Navier--Stokes--Fourier system is Levinson dissipative. Although the ballistic energy need be non--negative, we show that its
``entropy'' component is dominated by the internal energy. More specifically, the bound \eqref{w8} is equivalent to
\[
\limsup_{t \to \infty} \intO{ \left( \frac{1}{2} \vr |\vu - \vuB|^2 + \vr e \right) } \leq \mathcal{E}_\infty,
\]
modulo a suitable modification of $\mathcal{E}_\infty$.

Finally, note that if $\vtB$ is a positive \emph{constant}, then \eqref{w7} reduces to an energy inequality formally similar to that for the barotropic system studied in \cite{BreFeiNov20}. Thus, similarly to \cite{BreFeiNov20}, a hard sphere pressure equation of state is necessary to keep the density bounded and to control the first integral on the right--hand side of \eqref{w7}. On the other hand, under the non-slip boundary conditions
$\vuB = 0$, inequality \eqref{w7} gives rise to
	\begin{align}
	\frac{\D }{\dt} &\intO{ \left( \frac{1}{2} \vr |\vu|^2 + \vr e - \tvt \vr s \right) }   \br &+
	\intO{ \frac{\tvt}{\vt}	 \left( \mathbb{S} : \Ds \vu - \frac{\vc{q} \cdot \Grad \vt }{\vt} \right) }  \br
	&\leq
	-  \intO{ \left[ \vr s \left( \partial_t \tvt + \vu \cdot \Grad \tvt \right) + \frac{\vc{q}}{\vt} \cdot \Grad \tvt \right] };
	\label{w9}
\end{align}
whence the rightmost integral must be dominated by the dissipation
\[
\intO{ \frac{\tvt}{\vt}	 \left( \mathbb{S} : \Ds \vu - \frac{\vc{q} \cdot \Grad \vt }{\vt} \right) }.
\]
Again this does not seem realistic unless uniform bounds on the density $\vr$ are {\it a priori} imposed. The above arguments justify the choice of the hard sphere pressure equation of state introduced in the section below.

\section{Constitutive relations, equation of state}
\label{h}

Before stating the main results, we introduce the structural hypotheses imposed on the equations of state motivated by \cite{FeNo6A}. We suppose
\begin{equation} \label{h1}
	\widetilde{p}(\vr, \vt) = \vt^{\frac{5}{2}} P \left( \frac{\vr}{\vt^{\frac{3}{2}}  } \right) + \frac{a}{3} \vt^4,\
	\widetilde{e}(\vr, \vt) = \frac{3}{2} \frac{\vt^{\frac{5}{2}} }{\vr} P \left( \frac{\vr}{\vt^{\frac{3}{2}}  } \right) + \frac{a}{\vr} \vt^4, \ a > 0,
\end{equation}
where $P \in C^1[0,\infty)$ satisfies
\begin{equation} \label{h2}
	P(0) = 0,\ P'(Z) > 0 \ \mbox{for}\ Z \geq 0,\ 0 < \frac{ \frac{5}{3} P(Z) - P'(Z) Z }{Z} \leq c \ \mbox{for}\ Z > 0.
\end{equation} 	
In particular, the function $Z \mapsto P(Z)/ Z^{\frac{5}{3}}$ is decreasing, and we suppose
\begin{equation} \label{h3}
	\lim_{Z \to \infty} \frac{ P(Z) }{Z^{\frac{5}{3}}} = p_\infty \geq 0.
\end{equation}
The associated  entropy $s$ reads
\begin{equation} \label{h4}
	s(\vr, \vt) = \mathcal{S} \left( \frac{\vr}{\vt^{\frac{3}{2}} } \right) + \frac{4a}{3} \frac{\vt^3}{\vr},
\end{equation}
where
\begin{equation} \label{h5}
	\mathcal{S}'(Z) = -\frac{3}{2} \frac{ \frac{5}{3} P(Z) - P'(Z) Z }{Z^2}.
\end{equation}

In addition, following \cite{BreFeiNov20}, we introduce the hard sphere perturbation of the equation of state,
\begin{align}
	p(\vr, \vt) &= \widetilde{p}(\vr, \vt) + p_{\rm HS}(\vr), \br
	p_{\rm HS} &\in C^1[0, \Ov{\vr}),\ p_{\rm HS}(0) = 0,\ p_{\rm HS}' > 0 \ \mbox{in}\ (0, \Ov{\vr}),\
	\lim_{\vr \to \Ov{\vr}-} p_{\rm HS}(\vr) = \infty.
\label{h6}
	\end{align}
The related internal energy reads
\begin{equation} \label{h7}
	e(\vr, \vt) = \widetilde{e}(\vr, \vt) + \int_{\Ov{\vr}/2}^\vr \frac{p_{\rm HS}(z) }{z} \ \D z.
\end{equation}	

\begin{Remark} \label{Rh1}
	
	The main effect of the hard sphere pressure equation of state is, of course, the uniform bound imposed {\it a priori} on the fluid density,
	\begin{equation} \label{h7bis}
		0 \leq \vr(t,x) \leq \Ov{\vr},
		\end{equation}
for any weak solution of the Navier--Stokes--Fourier system.
	
	\end{Remark}

The transport coefficients $\mu$, $\eta$, and $\kappa$ are continuously differentiable functions of the temperature $\vt$ satisfying
\begin{align}
	0 < \underline{\mu} \left(1 + \vt^\Lambda \right) &\leq \mu(\vt) \leq \Ov{\mu} \left( 1 + \vt^\Lambda \right),\
	|\mu'(\vt)| \leq c \ \mbox{for all}\ \vt \geq 0,\ \frac{1}{2} \leq \Lambda \leq 1, \br
	0 &\leq  \eta(\vt) \leq \Ov{\eta} \left( 1 + \vt^\Lambda \right), \br
	0 < \underline{\kappa} \left(1 + \vt^\beta \right) &\leq  \kappa(\vt) \leq \Ov{\kappa} \left( 1 + \vt^\beta \right),\ \beta \geq 0.
	\label{h8}
\end{align}

In addition, we say that $s$ is compatible with the Third law of thermodynamics, if
\begin{equation} \label{h9}
	\lim_{\vt \to 0} s(\vr, \vt) = 0 \ \mbox{for any fixed}\ \vr > 0.
\end{equation}

The existence theory developed in \cite{ChauFei} can be easily modified to accommodate the hard--sphere pressure, at least in the specific form
\[
p_{\rm HS}(\vr) \approx (\Ov{\vr} - \vr)^{-\beta},\ \beta > 3,
\]
cf. also \cite{FeiZha}.

\section{Main results}
\label{m}

We are ready to present our main results. If not otherwise stated, we suppose that the boundary data enjoy the degree of smoothness necessary for the analysis, and that $\partial \Omega$ is sufficiently smooth.
In addition, without loss of generality, we suppose that the boundary data are restrictions of smooth functions defined on $(T, \infty) \times R^d$.

\subsection{General boundary conditions}

We suppose that $\Omega \subset R^d$, $d=2,3$ is a bounded domain of class $C^\infty$ such that
\begin{align}
	\partial \Omega = \cup_{i=0}^n \Gamma^i,\ \Gamma^i \cap \Gamma^j = \emptyset \ i \ne j,
	 \br
	\label{m1}
	\end{align}
where $\Gamma^i$ are connected components of $\partial \Omega$ and $\Gamma^0$ is the boundary of the unbounded component of $R^d \setminus \Omega$.

\begin{mdframed}[style=MyFrame]
	\begin{Theorem} \label{Tm1}
		Let $\Omega \subset R^d$, $d=2,3$ be a bounded domains of class $C^\infty$, the boundary of which admits the decomposition \eqref{m1}. In addition, suppose the following holds:
		
		\begin{itemize}
			\item The pressure $p$ and the internal energy $e$ are given by the hard sphere equations of state \eqref{h6}, \eqref{h7}, the entropy is compatible with the Third law of thermodynamics \eqref{h9}.
			
			\item The transport coefficients $\mu$, $\eta$, and $\kappa$ are continuously differentiable functions of $\vt$ satisfying \eqref{h8}, with $\Lambda = 1$, $\beta > 6$.
			
			\item The boundary data $\vrB$, $\vuB$, $\vtB$ are restrictions of continuously differentiable functions in
			$(T , \infty) \times R^d$, and
			\begin{align}
			0 &< \inf_{(T, \infty) \times R^d} \vrB \leq  \sup_{(T, \infty) \times R^d} \vrB \leq \Ov{\vr}, \br
			0 &< \underline{\vt} = \inf_{(T, \infty) \times R^d} \vtB \leq  \sup_{(T, \infty) \times R^d} \vtB = \Ov{\vt}, \br
			\int_{\Gamma^i} \vuB \cdot \vc{n} \ \D \sigma_x &= 0,\ i = 1,\dots,n,\
			\inf_{(T,\infty)} \int_{\Gamma_0} \vuB \cdot \vc{n} \ \D \sigma_x > 0, \br
			{ |\partial_t \vuB(t,x) | + |\partial_t \vtB(t,x)|} &+ | \Grad \vrB(t,x)| + |\Grad^{\alpha} \vuB(t,x) | + |\Grad^{\alpha} \vtB(t,x) | \leq D,\br
			\alpha &= 0,1,2,\ { t \in (T, \infty),\ x \in R^d.}
			\label{m2}
				\end{align}
			\item The driving force $\vc{g}$ is a bounded measurable function,
			\begin{equation} \label{m3}
			{\|\vc{g} \|_{L^\infty((T, \infty) \times \Omega; R^d)}  \leq D.}
			\end{equation}
			\end{itemize}
		
		Then there exists a universal constant $\mathcal{E}_\infty$, depending solely on the norm of the boundary data and the driving force, such that
		\begin{equation} \label{m3a}
		\limsup_{t \to \infty} \intO{ \left( \frac{1}{2} \vr |\vu - \vuB|^2 + \vr e (\vr, \vt) - \tvt \vr s (\vr, \vt) \right) } \leq \mathcal{E}_\infty
		\end{equation}
		for any weak solution $(\vr, \vt, \vu)$ of the Navier--Stokes--Fourier system on $(T, \infty) \times \Omega$,
	where $\tvt$ is the unique solutions of the Dirichlet problem
	\begin{equation} \label{m4}
		\Del \tvt(t, \cdot) = 0 \ \mbox{in}\ \Omega,\ \tvt (t, \cdot) |_{\partial \Omega} = \vtB (t, \cdot).
		\end{equation}
	\end{Theorem}
	\end{mdframed}

\medskip

\begin{Remark} \label{Rm1}
The condition
	\[
	 \int_{\Gamma_0} \vuB \cdot \vc{n} \ \D \sigma_x > 0
	\]
required for any $t \in (T, \infty)$ is purely ``compressible'' as it excludes the possibility of $\vuB$ being solenoidal. Its stabilizing effect has been observed in \cite{BreFeiNov20}.
It is worth noting that the same condition is imposed by Choe, Novotn\' y and Yang \cite[Theorem 2.5]{ChNoYa} to show the \emph{existence} of global in time weak solutions for the hard--sphere barotropic
Navier--Stokes system. This condition is relaxed in the forthcoming section.
	
	\end{Remark}

\bigskip

Note that the choice of parameters $\Lambda = 1$, $\beta > 6$ as well as \eqref{h9} are also necessary to show the \emph{existence} of global in time weak solutions in \cite[Theorem 4.2]{ChauFei}.
The bound \eqref{m3a} can be equivalently stated as
\begin{equation} \label{m3b}
	\limsup_{t \to \infty} \intO{ \left( \frac{1}{2} \vr |\vu - \vuB|^2 + \vr e (\vr, \vt)  \right) } \leq \mathcal{E}_\infty,
\end{equation}
without any reference to \eqref{m4}.

\subsection{No--slip boundary conditions, B\' enard problem}

Th reader will have noticed that the principal hypotesis \eqref{m2} of Theorem \ref{Tm1} does not include the no--slip boundary conditions $\vuB|_{\partial \Omega} = 0$. The next result focuses on the B\' enard problem, where the boundary temperature is prescribed, while the normal velocity vanishes on the boundary. In particular, the total mass
\[
M = \intO{ \vr }
\]
is a constant of motion.

\begin{mdframed}[style=MyFrame]
	\begin{Theorem} [Impermeable boundary] \label{Tm2}
Let $\Omega \subset R^d$, $d=2,3$ be a bounded domains of class $C^\infty$, the boundary of which admits the decomposition \eqref{m1}. In addition, suppose the following holds:

\begin{itemize}
	\item The pressure $p$ and the internal energy $e$ are given by the hard sphere equqations of state \eqref{h6}, \eqref{h7}, the entropy is compatible with the Third law of thermodynamics \eqref{h9}.
	
	\item The transport coefficients $\mu$, $\eta$, and $\kappa$ are continuously differentiable functions of $\vt$ satisfying \eqref{h8}, with $\Lambda = 1$, $\beta > 6$.
	
	\item The boundary data $\vuB$, $\vtB$ are restrictions of continuously differentiable functions in
	$(T , \infty) \times R^d$, and
	\begin{align}
		0 < \underline{\vt} &= \inf_{(T, \infty) \times R^d} \vtB \leq  \sup_{(T, \infty) \times R^d} \vtB = \Ov{\vt}, \br
		\vuB \cdot \vc{n}|_{\partial \Omega} &= 0, \br
		|\partial_t \vuB(t,x)| + |\partial_t \vtB(t,x)|   &+ |\Grad^{\alpha} \vuB(t,x) | + |\Grad^{\alpha} \vtB(t,x) | \leq D,\br
		\alpha &= 0,1,2,\  t \in (T, \infty), x \in R^d.	
		\label{m5}
	\end{align}
	\item The driving force $\vc{g}$ is a bounded measurable function,
	\begin{equation} \label{m6}
		\|\vc{g} \|_{L^\infty((T, \infty) \times \Omega)} \leq D.
	\end{equation}
\end{itemize}

Then there exists a universal constant $\mathcal{E}_\infty$, depending solely on the norm of the boundary data, the total mass $M$,
and the driving force, such that
\[
\limsup_{t \to \infty} \intO{ \left( \frac{1}{2} \vr |\vu - \vuB|^2 + \vr e_{HS} (\vr, \vt) - \tvt \vr s (\vr, \vt) \right) } \leq \mathcal{E}_\infty
\]
for any weak solution $(\vr, \vt, \vu)$ of the Navier--Stokes--Fourier system,
where $\tvt$ is the unique solutions of the Dirichlet problem \eqref{m4}.

\end{Theorem}		
\end{mdframed}

Note carefully that hypothesis \eqref{m5} is not a special case of \eqref{m2}. In particular, \eqref{m5} is compatible with the no-slip boundary condition $\vuB|_{\partial \Omega} = 0$.

\subsection{Convergence to equilibrium}

Finally, we discuss the situation, where the boundary temperature $\vtB$ is a positive constant, while the velocity field $\vuB$ coincides with a rigid motion tangential
to $\partial \Omega$. As $\Omega$ is bounded, this is possible only if:
\begin{itemize}
	\item $\vuB = 0$, or
	\item
$\vuB$ is a rigid rotation and $\Omega$ is radially symmetric with respect to the axis of rotation.
\end{itemize}

If the driving force $\vc{g} = \Grad G$ is potential, the global in time solutions are expected to converge to an equilibrium solution $(\vr_E, \vuB, \vtB)$, where
\begin{align}
	\Div (\vr_E \vuB) &= 0, \br
	\Div (\vr_E \vuB \otimes \vuB) + \Grad p(\vr_E, \vtB) &= \vr_E \Grad G, \br
	\intO{ \vr_E } = M,
\label{m7}	
\end{align}
where the total mass
\[
M = \intO{ \vr(t, \cdot) }
\]
is a constant of motion.

\begin{mdframed}[style=MyFrame]

	\begin{Theorem} [Convergence to equilibrium] \label{Tm3}
	Let $\Omega \subset R^d$, $d=2,3$ be a bounded domains of class $C^\infty$. In addition, suppose the following holds:
	
	\begin{itemize}
		\item The pressure $p = \widetilde{p}$ and the internal energy $e = \widetilde{e}$ are given by the constitutive equations \eqref{h1}--\eqref{h5}, with $p_\infty > 0$.
		
		\item The transport coefficients $\mu$, $\eta$, and $\kappa$ are continuously differentiable functions of $\vt$ satisfying \eqref{h8}, with $\frac{1}{2} \leq \Lambda \leq 1$, $\beta = 3$.
		
		\item The boundary data satisfy
		\begin{align}
		\vuB &= \vuB(x),\ \Ds \vuB = 0,\ \vuB \cdot \vc{n}|_{\partial \Omega} = 0, \br
		\vtB &> 0 \ \mbox{--a positive constant.}
		\label{m8}
		\end{align}
		\item The driving force $\vc{g}$ is potential,
		\begin{equation} \label{m9}
		\vc{g} = \Grad G, \ G= G(x),\ G \in W^{1,\infty}(\Omega),\ \Grad G \cdot \vuB = 0.
		\end{equation}
	\end{itemize}
	
	Then there exists a density profile $\vr_E$ solving the stationary problem \eqref{m7} such that
	\begin{align}
		\vr(t, \cdot) &\to \vr_E \ \mbox{in}\ L^{\frac{5}{3}}(\Omega), \br
		\vr \vu (t, \cdot) &\to \vr_E \vuB \ \mbox{in}\ L^{\frac{5}{4}}(\Omega; R^d), \br
		\vt(t, \cdot) &\to \vtB \ \mbox{in}\ L^4(\Omega)
		\end{align}
	as $t \to \infty$ for any weak solution $(\vr, \vt, \vu)$ of the Navier--Stokes--Fourier system.
\end{Theorem}		
\end{mdframed}

Unlike in Theorems \ref{Tm1}, \ref{Tm2}, the hard sphere pressure component is not necessary in Theorem \ref{Tm3}.

\section{Proof of the main results}
\label{u}

Our goal is to prove the main results stated in Section \ref{m}. First, as $\tvt$ is a solution of the Dirichlet problem \eqref{m4}, we may apply the standard maximum principle together with the elliptic estimates to deduce
\begin{align}
	0 < \underline{\vt} &\leq \tvt(t,x) \leq \Ov{\vt},\ \br
	\inf_{(T, \infty) \times
		\partial \Omega} \partial_t \vtB &\leq \partial_t \tvt (t,x) \leq \sup_{(T, \infty) \times
		\partial \Omega} \partial_t \vtB,\br
	|\Grad \tvt (t,x) | &\leq c(D) \ \mbox{for all}\ t > T,\ x \in \Omega.
	\label{u1}
	\end{align}

Next, we adapt the construction of a suitable extension of the velocity field $\vuB$ used in \cite{BreFeiNov20} to the time--dependent setting.
It is easy to observe that the component $\Gamma^0$ of the boundary $\partial \Omega$ contains at least one extremal point
$x_0 \in \Gamma_0$ satisfying
\[
\Ov{\Omega} \cap \tau_{x_0} = x_0,\ \mbox{where}\ \tau_{x_0}
\ \mbox{denotes the tangent plane to}\ \partial \Omega \ \mbox{at}\ x_0.
\]
Without loss of generality, we may assume that
\[ \Omega \subseteq \{x: x^1 <x_0^1\} \mbox{ and }
x_0 = [x_0^1,0,\dots,0],\ \tau_{x_0} = x_0 + R^{d-1}.
\]
Now,
consider a function
\[
\chi(z) = \left\{ \begin{array}{l} \ 0 \ \mbox{if}\ z \leq 0,\\
	\chi'(z) > 0 \ \mbox{for}\ z > 0 \end{array} \right.,
\]
together with a vector field
\[
\vc{v}^0_B(t,x) = \lambda(t) \left[ \chi (x^1 - x_0^1 + \delta), 0 , \dots, 0 \right].
\]
It is easy to check that
\[
\Ds \vc{v}^0_B = \begin{bmatrix} \lambda (t) \chi' (x^1 - x_0^1 + \delta)  & 0 & 0 \\
	0 & 0 &0 \\
	0 & 0 & 0
\end{bmatrix},\ \Div \vc{v}^0_B = \lambda (t) \chi' (x^1 - x_0^1 + \delta).
\]
Next, choose $\delta > 0$ small enough so that
\[
\vc{v}^0_B |_{\Gamma_i } = 0, \ i =1, \dots, n,\ \vc{v}^0_b|_{\Gamma _0} \ne 0,
\]
and then $\lambda (t) > 0$ large enough so that
\[
\int_{\Gamma^0} \vc{v}^0_B \cdot \vc{n} \ \D \sigma_x = \int_{\Gamma^0} \vuB \cdot \vc{n}\ \sigma_x > 0
\]
in accordance with hypothesis \eqref{m2}. We decompose
\[
\vuB = \vc{w}_B + \vc{v}^0_B,
\]
where
\begin{align}
	\int_{\Gamma^i} \vc{w}_B \cdot \vc{n} \ \D \sigma_x &= 0 \ \mbox{for}\ i=0,1,\dots, n, \br
	\Ds \vc{v}^0_B &\geq 0,\
	\ \mbox{there is an open set}\ B \subset \Omega, \ |B| > 0,\
	\inf_{(T, \infty) \times B} (\Div \vc{v}^0_B) \geq \underline{d} > 0.
	\label{u2}
	\end{align}
Under the hypotheses of Theorem \ref{Tm2}, where, obviously,
\[
\int_{\Gamma^0} \vuB \cdot \vc{n} \ \D \sigma_x = 0,
\]	
we simply set $\vc{v}^0_B = 0$, $\vc{w}_B = \vuB$.

Now, exactly as in \cite[Section 4.2]{BreFeiNov20}, we use Galdi \cite[Lemma IX.4.1]{GALN},
Kozono and Yanagisawa \cite[Proposition 1]{KozYan} to write
\begin{equation} \label{u3}
	\vc{w}_B (t,x) = {\bf curl}_x (d_\ep (x) \vc{z}_B(t,x)),\ {\bf curl_x}\  \vc{z}_B = \vc{w}_B,
	\end{equation}
where
\begin{align}
	|d_\ep | &\leq 1,\ d_\ep (x) \equiv 1 \ \mbox{for all}\ x \ \mbox{in an open neighborhood of}\ \partial \Omega, \br
	d_\ep (x) &\equiv 0 \ \mbox{whenever}\ {\rm dist}[x, \partial \Omega] > \ep, \br
	|D^\alpha_x d_\ep (x) | &\leq c \frac{\ep}{{\rm dist}^{|\alpha|} [x, \partial \Omega]},\ |\alpha| = 1,2 ,\ 0 < \ep < 1,\
	x \in \Omega.
\label{u4}	
\end{align}

\begin{Remark} \label{mR1}
	If $d=2$,
	we adopt the convention that $z_B$ is scalar and the operator ${\bf curl}_x$ is replaced by $\Grad^\perp$ in \eqref{u3}.
	\end{Remark}

\subsection{Proof of Theorems \ref{Tm1} \ref{Tm2}}

First, write the ballistic energy inequality \eqref{w6} in the form
	\begin{align}
	&\left[ \intO{ \left( \frac{1}{2} \vr |\vu - \vuB|^2 + \vr e - \tvt \vr s \right) } \right]_{t = I}^{t = I + \tau}  \br &+
	\int_I^{I + \tau} \intgi{ \Big[ \vrB e (\vrB, \vtB) - \vtB \vrB s(\vrB, \vtB) \Big] } \dt \br &+
	\int_I^{I + \tau} \intgo{ \Big[ \vr e(\vr, \vtB) - \vtB \vr s(\vr, \vtB) \Big] } \dt \br &+
	\int_I^{I + \tau}\intO{ \frac{\tvt}{\vt}	 \left( \mathbb{S}(\vt, \Ds \vu) : \Ds \vu - \frac{\vc{q} \cdot \Grad \vt }{\vt} \right) } \dt \br
	\leq
	&-  \int_I^{I + \tau} \intO{ \Big[ \vr (\vu - \vuB) \otimes (\vu - \vuB) + p \mathbb{I} - \mathbb{S} \Big] : \Ds \vuB } \dt \br &+
	\int_I^{I + \tau} \intO{ \vr (\vu - \vuB)\cdot (\vc{g} - \partial_t \vuB - \vuB \cdot \Grad \vuB) } \dt \br
	&- \int_I^{I + \tau} \intO{ \left[ \vr s \left( \partial_t \tvt + \vu \cdot \Grad \tvt \right) + \frac{\vc{q}}{\vt} \cdot \Grad \tvt \right] } \dt,\ I > T, \ \tau \geq 0.
	\label{u5}
\end{align}

As $\Lambda = 1$ in hypothesis \eqref{h8}, we may use Korn--Poincar\' e inequality to estimate
\begin{align}
\| \vu - \vuB \|_{W^{1,2}_0 (\Omega; R^d)}^2 &\aleq \left\| \Grad (\vu - \vuB) + \Grad^t (\vu - \vuB) - \frac{2}{d} \Div (\vu - \vuB)\mathbb{I} \right\|_{L^{2}(\Omega; R^{d \times d})}^2, \br
&\aleq \intO{ \left| \Grad \vu  + \Grad^t \vu  - \frac{2}{d} \Div \vu \mathbb{I} \right|^2 } + c(\| \Grad \vuB \|_{L^2}), \br
&\aleq \intO{ \frac{\tvt}{\vt} \mathbb{S}(\vt, \Ds \vu ): \Ds \vu }   + c(\| \Grad \vuB \|_{L^2}),
\nonumber
\end{align}
where $A\aleq B$ means $A\leq c B$ for a constant number $c\geq0$.
Consequently, as $\vr$ is bounded by $\Ov{\vr}$ uniformly in \eqref{h7bis},
	\begin{align}
	&\int_I^{I + \tau} \intO{ \vr (\vu - \vuB)\cdot (\vc{g} - \partial_t \vuB - \vuB \cdot \Grad \vuB) } \dt\nonumber\\
	&\qquad\qquad\leq 	\ep \int_I^{I + \tau}  \left\| \vu - \vuB \right\|^2_{W^{1,2}_0(\Omega; \R^d)}\dt+c \Big(\tau, \ep ,D,  \Ov{\vr}, \| \vc{g} \|_{L^\infty} \Big)
	\label{u5-1}
\end{align}
for any $\ep > 0$. Thus we deduce from
from \eqref{u5}, \eqref{u5-1}
	\begin{align}
	&\left[ \intO{ \left( \frac{1}{2} \vr |\vu - \vuB|^2 + \vr e - \tvt \vr s \right) } \right]_{t = I}^{I = I + \tau}  \br  &+
	\int_I^{I + \tau} \left\| \vu - \vuB \right\|^2_{W^{1,2}_0(\Omega; R^d)} \dt + \int_{I}^{I + \tau} \intO{ \frac{\kappa (\vt) |\Grad \vt|^2 }{\vt^2}  } \dt \br
	\aleq
	&-  \int_I^{I + \tau} \intO{ \Big[ \vr (\vu - \vuB) \otimes (\vu - \vuB) + p \mathbb{I} - \mathbb{S} \Big] : \Ds \vuB } \dt \br
	&- \int_I^{I + \tau} \intO{ \left[ \vr s \left( \partial_t \tvt + \vu \cdot \Grad \tvt \right) - \frac{\kappa (\vt) \Grad \vt }{\vt} \cdot \Grad \tvt \right] } + c \Big(\tau, D, \underline{\vr}, \Ov{\vr}, \underline{\vt}, \Ov{\vt}, \| \vc{g} \|_{L^\infty} \Big).
	\label{u6}
\end{align}

Next, as $\tvt$ solves \eqref{m4}, the Gauss--Green integration formula yields
\begin{equation} \label{u7}
\intO{ \frac{\kappa (\vt) \Grad \vt }{\vt} \cdot \Grad \tvt } = \intO{ \Grad \mathcal{K}(\vt) \cdot \Grad \tvt } = \int_{\partial \Omega} \mathcal{K} (\vtB) \Grad \tvt \cdot \vc{n} \ \D \sigma_x,
\end{equation}
where $\mathcal{K}'(\vt) = \frac{\kappa (\vt) }{\vt}$. Furthermore,
\begin{equation} \label{u8}
\intO{ \mathbb{S} (\vt, \Ds \vu ) : \Ds \vuB } \aleq \ep \intO{ |\Grad \vu |^2 |\Ds \vuB| }  + c(\ep) \intO{ (1 + \vt^2) |\Ds \vuB| }
\end{equation}
for any $\ep > 0$. Finally, as the entropy is compatible with the Third law of thermodynamics \eqref{h9}, we have
\begin{equation} \label{u9}
	|\vr s(\vr, \vt) | \aleq \left( 1 + \vr |\log (\vr) | + \vr [\log(\vt)]^+ + \vt^3 \right).
	\end{equation}
In view of \eqref{u7}--\eqref{u9}, the inequality \eqref{u6} gives rise to
	\begin{align}
	&\left[ \intO{ \left( \frac{1}{2} \vr |\vu - \vuB|^2 + \vr e - \tvt \vr s \right) } \right]_{t = I}^{t = I + \tau}  \br  &+
	\int_I^{I + \tau} \left\| \vu - \vuB \right\|^2_{W^{1,2}_0(\Omega; R^d)} \dt + \int_{I}^{I + \tau} \intO{ \frac{\kappa (\vt) |\Grad \vt|^2 }{\vt^2}  } \dt \br
	\aleq
	&-  \int_I^{I + \tau} \intO{ \Big[ \vr (\vu - \vuB) \otimes (\vu - \vuB) \Big] :  \Ds \vuB } \dt  - \int_I^{I + \tau} \intO{ p \Div \vuB } \dt \br
	&+ \int_I^{I + \tau} \intO{ \vt^3 (1 + |\vu| ) } \dt + c \Big(\tau, D, \underline{\vr}, \Ov{\vr}, \underline{\vt}, \Ov{\vt}, \| \vc{g} \|_{L^\infty} \Big).
	\label{u10}
\end{align}

Now, in accordance with hypothesis \eqref{h8},
\[
	\| \Grad \log (\vt) \|^2_{L^2(\Omega; R^d)} + \| \Grad \vt^{\frac{\beta}{2}} \|^2_{L^2(\Omega; R^d)} \aleq \intO{ \frac{\kappa (\vt) |\Grad \vt|^2 }{\vt^2}  }.
\]
Thus, as the boundary values of $\vt$ are controlled,
\begin{equation} \label{u11}
\| \log (\vt) \|^2_{W^{1,2}(\Omega)} + \| \vt^{\frac{\beta}{2}} \|^2_{W^{1,2}(\Omega)} \aleq \intO{ \frac{\kappa (\vt) |\Grad \vt|^2 }{\vt^2}  } + c( \underline{\vt}, \Ov{\vt})
\end{equation}

Next, we have
\[
\intO{ \vt^3 |\vu| } \leq \ep \| \vu \|^2_{L^2(\Omega; R^d)} + c(\ep) \| \vt^3 \|^2_{L^2(\Omega)}
\]
and, since $\beta > 6$,
\begin{equation} \label{u12}
\| \vt^3 \|^2_{L^2(\Omega)} \leq \delta \| \vt^{\frac{\beta}{2}} \|_{L^2(\Omega)}^2 + c(\delta)
\end{equation}
for any $\ep > 0$, $\delta > 0$.
Consequently, going back to \eqref{u10} we may infer that
\begin{align}
	&\left[ \intO{ \left( \frac{1}{2} \vr |\vu - \vuB|^2 + \vr e - \tvt \vr s \right) } \right]_{t = I}^{t = I + \tau}  \br  &+
	\int_I^{I + \tau} \left\| \vu - \vuB \right\|^2_{W^{1,2}_0(\Omega; R^d)} \dt + \int_{I}^{I + \tau}\left(  \| \log (\vt) \|^2_{W^{1,2}(\Omega)} + \| \vt^{\frac{\beta}{2}} \|^2_{W^{1,2}(\Omega)} \right)  \dt \br
	& + \int_I^{I + \tau} \intO{ p \Div \vc{v}_B^0 } \dt \br
	\aleq
	&-  \int_I^{I + \tau} \intO{ \Big[ \vr (\vu - \vuB) \otimes (\vu - \vuB) \Big] :  \Ds \vuB } \dt + c \Big(\tau, D, \underline{\vr}, \Ov{\vr}, \underline{\vt}, \Ov{\vt}, \| \vc{g} \|_{L^\infty} \Big).
	\label{u13}
\end{align}

Finally, exactly as in \cite[Section 5]{BreFeiNov20}, it can be shown that the integral
\[
\intO{ \Big[ \vr (\vu - \vuB) \otimes (\vu - \vuB) \Big] :  \Ds \vuB }
\]
can be absorbed by the left--hand side of \eqref{u13} thanks to the decomposition \eqref{u2}, \eqref{u3} as long as $\ep > 0$ in \eqref{u3} is fixed small enough. Thus, by virtue of \eqref{u2}, we conclude
\begin{align}
	&\left[ \intO{ \left( \frac{1}{2} \vr |\vu - \vuB|^2 + \vr e - \tvt \vr s \right) } \right]_{t = I}^{t = I + \tau}  \br  &+
	\int_I^{I + \tau} \left\| \vu - \vuB \right\|^2_{W^{1,2}_0(\Omega; R^d)} \dt + \int_{I}^{I + \tau}\left(  \| \log (\vt) \|^2_{W^{1,2}(\Omega)} + \| \vt^{\frac{\beta}{2}} \|^2_{W^{1,2}(\Omega)} \right)  \dt \br
	& + \underline{d} \int_I^{I + \tau} \int_B p \ \dx \dt
	\leq c \Big(\tau, D, \underline{\vr}, \Ov{\vr}, \underline{\vt}, \Ov{\vt}, \| \vc{g} \|_{L^\infty} \Big).
	\label{u14}
\end{align}
Note that the same inequality  with $\underline{d} = 0$ is obtained under the hypotheses of Theorem \ref{Tm2}. Our ultimate goal is therefore to derive the estimate
\begin{equation} \label{u15}
\intO{ \vr e } \aleq 	\left\| \vu - \vuB \right\|^2_{W^{1,2}_0(\Omega; R^d)} + \| \log (\vt) \|^2_{W^{1,2}(\Omega)} + \| \vt^{\frac{\beta}{2}} \|^2_{W^{1,2}(\Omega)} +
\int_B p(\vr, \vt) \dx.
\end{equation}	
Of course, the main problem is that the last integral is evaluated only over the ``small'' set $B$ identified in \eqref{u2}.

\subsubsection{Pressure estimates, proof of Theorem \ref{Tm1}}

Similarly to \cite{BreFeiNov20},
the desired bound \eqref{u15} follows from the pressure estimates. As the present setting is slightly different due to the temperature depending terms, we reproduce some details of the proof for reader's convenience.
First, we recall the Bogovskii operator:
\begin{align}
	\mathcal{B} : L^q_0 (\Omega) &\equiv \left\{ f \in L^q(\Omega) \ \Big| \ \intO{ f } = 0 \right\}
	\to W^{1,q}_0(\Omega, R^d),\ 1 < q < \infty,\br
	\Div \mathcal{B}[f] &= f ,
	\nonumber
\end{align}
see e.g. Galdi \cite[Chapter 3]{GALN}, Geissert, Heck, and Hieber \cite{GEHEHI}.

Using $\psi(t) \mathcal{B}[\Phi]$, $\Phi = \Phi(x)$ as a test function in the momentum equation \eqref{w4}, we deduce
\begin{align}
		\int_I^{I + \tau} &\intO{ p (\vr, \vt) \Phi } \dt =
		\left[ \intO{ \vr \vu \cdot \mathcal{B} \left[ \Phi \right] }  \right]_{t=I}^{t = I + \tau} \br
		&- \int_I^{I + \tau} \intO{  \vr \vu \otimes \vu : \Grad \mathcal{B} \left[ \Phi \right] } \dt
		+ \int_I^{I + \tau} \intO{  \mathbb{S} ( \Ds \vu ) : \Grad \mathcal{B} \left[\Phi \right] } \dt \br
		&- \int_I^{I + \tau} \intO{  \vr \vc{g} \cdot \mathcal{B} \left[ \Phi \right] } \dt.
	\label{u16}
\end{align}
All integrals on the right--hand side can be estimated exactly as in \cite[Section 5]{BreFeiNov20} except
\[
\intO{ \mathbb{S} (\vt, \Ds \vu ) : \Grad \mathcal{B} \left[\Phi \right]  }
\]
that can be treated as
\begin{align}
\left| \intO{ \mathbb{S} (\vt, \Ds \vu ) : \Grad \mathcal{B} \left[\Phi \right]  } \right| &\aleq \|  \Grad \mathcal{B} \left[\Phi \right] \|_{L^4(\Omega)} \left(
\| \Grad \vu \|^2_{L^2(\Omega; R^{d \times d})} + \| (1 + \vt) \|_{L^4(\Omega)}^2 \right) \br
\aleq  \|  \Grad \mathcal{B} \left[\Phi \right] \|_{L^4(\Omega)} &\left( \left\| \vu - \vuB \right\|^2_{W^{1,2}_0(\Omega; R^d)} + \| \vt^{\frac{\beta}{2}} \|^2_{W^{1,2}(\Omega)}
+ c \Big(D, \underline{\vr}, \Ov{\vr}, \underline{\vt}, \Ov{\vt}\Big) \right).
\label{u17}
\end{align}

With the estimate \eqref{u17} at hand, we can repeat step by step the arguments of \cite[Section 5]{BreFeiNov20} to deduce from \eqref{u14}, \eqref{u15} the final bound
\begin{align}
	&\left[ \intO{ \left( \frac{1}{2} \vr |\vu - \vuB|^2 + \vr e - \tvt \vr s -  \vr \vu \cdot \mathcal{B}[\Phi] \right) } \right]_{t = I}^{t = I + \tau}  \br  &+
	\delta \int_I^{I + \tau} \left\| \vu - \vuB \right\|^2_{W^{1,2}_0(\Omega; R^d)} \dt + \delta \int_{I}^{I + \tau}\left(  \| \log (\vt) \|^2_{W^{1,2}(\Omega)} + \| \vt^{\frac{\beta}{2}} \|^2_{W^{1,2}(\Omega)} \right)  \dt \br
	& +  \delta \int_I^{I + \tau} \intO{ p } \dt
	\leq c \Big(\tau, D, \underline{\vr}, \Ov{\vr}, \underline{\vt}, \Ov{\vt}, \| \vc{g} \|_{L^\infty}, \Phi \Big),
	\label{u18}
\end{align}
for a suitable $\Phi$ and $\delta > 0$.

Similarly to \cite[Section 5]{BreFeiNov20}, relation \eqref{u18} yields the conclusion of Theorem \ref{Tm1}. First observe that
	\begin{equation} \label{u18ab}
	\vr e \aleq p(\vr, \vt)
	\end{equation}
	see \cite[formula (2.4)]{BreFeiNov20}.
	
	Denote
\begin{equation} \label{u18b}
\mathcal{E}(t) = \intO{ \left( \frac{1}{2} \vr |\vu - \vuB|^2 + \vr e - \tvt \vr s -  \vr \vu \cdot \mathcal{B}[\Phi] \right) (t, \cdot) }
\end{equation}		
As the entropy $s$ satisfies \eqref{u9},
\begin{equation} \label{u18bb}
\mathcal{E} \geq \lambda  \intO{ \left( \frac{1}{2} \vr |\vu - \vuB|^2 + \vr e \right) } - \frac{1}{\lambda}
\end{equation}
for some $\lambda > 0$.
Consequently, inequality \eqref{u18} yields
\begin{equation} \label{uu20}
[\mathcal{E}(t) ]_{t = I}^{t= I + \tau} + \delta \int_{I}^{I + \tau} \mathcal{E}(s) \ \D s \leq C(\tau, \delta) \ \mbox{for some}\ \delta > 0.
\end{equation}
For $\tau = 1$ we get the following dichotomy:
\begin{enumerate}
	\item Either there exists $t \in [I, I+ 1]$ such that
	\begin{equation} \label{uu18}
	\mathcal{E}(t) \leq \frac{2 C(1, \delta)}{\delta};
\end{equation}
\item
or
\begin{equation} \label{uu19}
\mathcal{E}(I + 1) \leq \mathcal{E}(I) - C(1, \delta).
\end{equation}
\end{enumerate}
In virtue of \eqref{u18bb},  $\mathcal{E}$ is bounded below and thus it follows that there exists $x_t\in [t-1,t]\subset (T,\infty)$ such that
\[
\mathcal{E}(x_t) \leq \frac{2C(1, \delta)}{\delta}.
\]
From \eqref{uu20},
\[
\mathcal{E}(t)\leq \mathcal{E}(x_t)+C(1,\delta)\leq \frac{2C(1, \delta)}{\delta}+C(1,\delta).
\]
So we get
\[
\limsup_{t\rightarrow\infty}\mathcal{E}(t)\leq C(1,\delta)\Big(1+\frac{2}{\delta}\Big).
\]
Consequently, we prove that
\[
\limsup_{t\rightarrow\infty} \intO{ \left( \frac{1}{2} \vr |\vu - \vuB|^2 + \vr e \right) } \leq \frac{C(1,\delta)}{\lambda}\Big(1+\frac{2}{\delta}\Big)+\frac{1}{\lambda^2}:=\mathcal{E}_\infty,
\]
which completes to prove Theorem \ref{m1}.

\subsubsection{Pressure estimates, proof of Theorem \ref{Tm2}}

Under the hypotheses of Theorem \ref{Tm2}, the total mass of the fluid is conserved,
\[
\intO{ \vr(t, \cdot) } = M \ \mbox{for any}\ t.
\]
Keeping in mind \eqref{u18}, the proof of Theorem \ref{Tm2} can be done by the same arguments as in \cite[Section 6.1]{BreFeiNov20}. Specifically, we repeat the pressure estimates with the test function
\[
\bfphi = \psi(t) \mathcal{B} \left[ \vr - \frac{1}{|\Omega|} \intO{ \vr } \right].
\]

\begin{mdframed}[style=MyFrame]
	
	\begin{Remark} \label{Ru1}
	It is worth noting that the above proof does not use any \emph{structural properties} of $p$ and $e$ as soon as the uniform bound on the density is established. The only piece of information to be retained being $\vr e \approx p$. In particular, the radiation component is irrelevant and possibly a more realistic equation of state similar to \cite{CiFeJaPe1} can be used.
	
	\end{Remark}
	
	\end{mdframed}

\subsection{Proof of Theorem \ref{Tm3}}

First we claim that for a given velocity field $\vuB$, the temperature $\vtB$, and the total mass $M > 0$, there exists a unique density profile $\vr_E$ solving the stationary problem \eqref{m7}.
Indeed, as $\Ds \vuB = 0$, the convective term in \eqref{m7} reads
\[
\Div (\vr_E \vuB \otimes \vuB) = \Div (\vr_E \vuB) \vuB +
\vr_E \vuB \cdot \Grad \vu_B = - \frac{1}{2} \vr_E \Grad |\vuB|^2.
\]
Accordingly, the problem \eqref{m7} can be rewritten as
\begin{align}
	\Grad \vr_E \cdot \vuB &= 0, \br
	\Grad p(\vr_E, \vtB) &= \vr_E \left( \Grad G + \frac{1}{2} |\vuB |^2 \right), \br
	\intO{ \vr_E } &= M.
	\label{u18a}	
\end{align}
As $P'(0) > 0$ in \eqref{h2} and $\vtB$ is given, the desired uniqueness result follows from
\cite[Theorem 2.1]{FP7}.

Under the hypotheses of Theorem \ref{Tm3}, the ballistic energy inequality \eqref{w7} simplifies considerably:
	\begin{align}
	\frac{\D }{\dt} &\intO{ \left( \frac{1}{2} \vr |\vu - \vuB|^2 + \vr e - \vtB \vr s \right) }   \br  &+
	\intO{ \frac{\tvt}{\vt}	 \left( \mathbb{S} : \Ds \vu - \frac{\vc{q} \cdot \Grad \vt }{\vt} \right) }  \br
	&\leq
	\intO{ \vr (\vu - \vuB)\cdot (\Grad G - \vuB \cdot \Grad \vuB) }  .
	\label{u19}
\end{align}
Moreover, by virtue of hypothesis \eqref{m9},
\[
	\intO{ \vr (\vu - \vuB)\cdot (\Grad G - \vuB \cdot \Grad \vuB) } = \intO{ \vr \vu \cdot \left(\Grad G + \frac{1}{2} \Grad |\vuB|^2 \right) };
\]
whence \eqref{u19} reduces to
	\begin{align}
	\frac{\D }{\dt} &\intO{ \left( \frac{1}{2} \vr |\vu - \vuB|^2 + \vr e - \vr \left( G + \frac{1}{2}  |\vuB|^2 \right) - \vtB \vr s \right) }   \br &+
	\intO{ \frac{\vtB}{\vt}	 \left( \mathbb{S}(\vt, \Ds \vu) : \Ds \vu + \frac{\kappa (\vt) |\Grad \vt|^2 }{\vt} \right) }  \
	\leq  0.
	\label{u20}
\end{align}

It turn out that the modified ballistic energy
\[
\intO{ \left( \frac{1}{2} \vr |\vu - \vuB|^2 + \vr e - \vr \left( G + \frac{1}{2}  |\vuB|^2 \right) - \vtB \vr s \right) }
\]
is a Lyapunov function decreasing along trajectories for which
\[
\intO{ \frac{\vtB}{\vt}	 \left( \mathbb{S}(\vt, \Ds \vu) : \Ds \vu + \frac{\kappa (\vt) |\Grad \vt|^2 }{\vt} \right) } > 0.
\]
In particular, there holds
\begin{align}
\int_T^\infty &\intO{ \frac{\vtB}{\vt}	 \left( \mathbb{S}(\vt, \Ds \vu) : \Ds \vu + \frac{\kappa (\vt) |\Grad \vt|^2 }{\vt} \right) } \dt \ \br
&= \int_T^\infty \intO{ \frac{\vtB}{\vt}	 \left( \mathbb{S}(\vt, \Ds \vu - \Ds \vuB) : (\Ds \vu - \Ds \vuB) + \frac{\kappa (\vt) |\Grad \vt - \Grad \vtB|^2 }{\vt} \right) } \dt < \infty
\label{u22}
\end{align}
for any weak solution of the Navier--Stokes--Fourier system defined on the time interval $(T, \infty)$.

Let $T_n \to \infty$ be a sequence of time. Let
\[
\vr_n (t,x) = \vr(T_n + t, x),\ \vt_n(t,x) = \vt(T_n + t,x),\ \vu_n(t,x) = \vu(T_n + t,x)
\]
be the associated time shifts of a global in time weak solution to the Navier--Stokes--Fourier
system. It follows from \eqref{u20}, \eqref{u22} that
\begin{align}
\vt_n &\to \vtB \ \mbox{in}\ L^{4 + \alpha}((0, T) \times \Omega),\br
\vu_n &\to \vuB \ \mbox{in}\ L^\alpha (0,T; W^{1,\alpha}(\Omega; R^2))
\nonumber
\end{align}
as $n \to \infty$ for some $\alpha > 1$. In particular,
\[
\Div \vu_n \to 0 \ \mbox{in}\ L^\infty((0,T) \times \Omega),
\]
which yields (cf. \cite[Chapter 4, Theorem 4.2]{FeiPr})
\begin{equation} \label{u23}
\vr_n \to \vr \ \mbox{in} \ L^{\frac{5}{3} + \alpha}((0,T) \times \Omega),
\end{equation}
passing to a suitable subsequence as the case may be.

Our final claim is that $\vr = \vr_E$, in particular, the convergence in \eqref{u23} is unconditional, which completes the proof of Theorem \ref{Tm3}. Seeing that the limit is again a weak solution of the Navier--Stokes--Fourier system we get
\begin{align}
\partial_t \vr + \Div (\vr \vuB) &= 0, \br
\partial_t (\vr \vuB) + \Div (\vr \vuB \otimes \vuB) + \Grad p(\vr, \vtB) &= \vr \Grad G
 \nonumber
\end{align}
in $\mathcal{D}'((0,T) \times \Omega)$.
Consequently,
\[
\Grad p(\vr, \vtB) = \vr \Grad \left( G + \frac{1}{2} |\vuB|^2 \right) \ \mbox{for any}\
t \in (0,T).
\]
Since $\intO{ \vr(t, \cdot) } = M$, the uniqueness result \cite{FP7} yields
$\vr(t, \cdot) = \vr_E$ for a.a. $t \in (0,T)$.

\section{Concluding remarks}
\label{a}

The hypotheses of Theorems \ref{Tm1}, \ref{Tm2} can be slightly relaxed. In the presence of the
singular hard--sphere pressure component, the structural hypotheses \eqref{h1}--\eqref{h5}
are not necessary. We may consider $p$, $e$ in the form \eqref{h6}, \eqref{h7}, where
$p$, $e$ are related to the entropy $s$ through general Gibbs' equation
\[
\vartheta D s = D e + p D \left( \frac{1}{\vr} \right)
\]
and satisfy the hypothesis of thermodynamic stability
\[
\frac{\partial p}{\partial \vr} > 0,\ \frac{\partial e}{\partial \vt} > 0.
\]
Strictly speaking, the presence of the radiation pressure is not necessary for the results of this paper, however, it is essential for the \emph{existence} of weak solutions.

The structural hypotheses that guarantee the existence of bounded absorbing sets are expected to provide positive results concerning the qualitative behavior of solutions in the long run.
In particular:

\begin{itemize}
	\item The existence of \emph{time periodic solutions} for problems driven by time periodic
	boundary data.
	\item The existence of \emph{global attractors}.
	\item Convergence of the ergodic averages and the existence of \emph{statistical stationary solutions} in the spirit of Constantin and Wu \cite{ConWu}, Foias, Rosa, Temam
	\cite{FoRoTe2}, \cite{FoRoTe1}, Vishik and Fursikov \cite{VisFur}.
	\end{itemize}

\noindent
These issues will be addressed in the future work.

\def\cprime{$'$} \def\ocirc#1{\ifmmode\setbox0=\hbox{$#1$}\dimen0=\ht0
	\advance\dimen0 by1pt\rlap{\hbox to\wd0{\hss\raise\dimen0
			\hbox{\hskip.2em$\scriptscriptstyle\circ$}\hss}}#1\else {\accent"17 #1}\fi}


\end{document}